\documentclass[11pt]{amsart}
\usepackage{fullpage}
\usepackage{graphicx}
\usepackage{amssymb}
\usepackage{epstopdf}
\usepackage{comment}
\newtheorem{theorem}{Theorem}[section]
\newtheorem{corollary}[theorem]{Corollary}
\newtheorem{proposition}[theorem]{Proposition}

\newtheorem{lemma}[theorem]{Lemma}
\newtheorem{problem}[theorem]{Problem}

\newtheorem{claim}[theorem]{Claim}

\numberwithin{equation}{section}

\def\S{\mathcal{S}}

\def\d{\delta}
\def\D{\Delta}

\def\supp{\text{supp}}

\def\COMMENT#1{}
\def\al#1{}
	\renewcommand{\al}[1]{\COMMENT{\textbf{AL: }#1}}  

\newenvironment{proofclaim}[1][Proof of claim]{\begin{proof}[#1]}{\end{proof}}

\title{Complete subgraphs in a multipartite graph}
\author{Allan Lo, Andrew Treglown, Yi Zhao}
\thanks{AL: University of Birmingham, United Kingdom, {\tt s.a.lo@bham.ac.uk}, research supported by EPSRC grant EP/V002279/1. \\
\indent AT: University of Birmingham, United Kingdom, {\tt a.c.treglown@bham.ac.uk}, research supported by EPSRC grant EP/V002279/1.\\
\indent YZ: Georgia State University, USA, {\tt yzhao6@gsu.edu}, research supported by NSF grant DMS 1700622 and Simons Collaboration Grant 710094.}
      
\begin{document}
\begin{abstract}
In 1975 Bollob\'as, Erd\H os, and Szemer\'edi asked the following question: given positive integers $n, t, r$ with $2\le t\le r-1$, what is the largest minimum degree $\delta(G)$ among all $r$-partite graphs $G$ with parts of size $n$ and which do not contain a copy of $K_{t+1}$? The $r=t+1$ case has attracted a lot of attention and was fully resolved by Haxell and Szab\'{o}, and Szab\'{o} and Tardos in 2006. In this paper we investigate the $r>t+1$ case of the problem, which has remained dormant for over forty years. We resolve the problem exactly in the case when $r \equiv -1 \pmod{t}$, and up to an additive constant for  many other cases, including when $r \geq (3t-1)(t-1)$.
%\footnote{YZ removed a sentence because it is not very strong: \emph{In all remaining open cases of the problem we provide constructions that improve on the previously best known extremal constructions}.}
Our approach utilizes a connection to the related problem of determining the maximum of the minimum degrees among the family of balanced $r$-partite $rn$-vertex graphs of chromatic number at most $t$.
\end{abstract}

\maketitle

\section{Introduction}
%In 1972 Erd\H os asked the $r=t+1$ case of the following question~\cite[Problem 2, 353--354]{MR0329905}.
\COMMENT{{\bf AT: rewritten first sentence as arguably Mantel's theorem started extremal gt, but Turan's theorem is the result that `underpins' the topic}}
The foundation stone of extremal graph theory is Tur\'an's theorem from 1941~\cite{Tu41}, which states that the Tur\'an graph $T_t(n)$ (the complete $t$-partite graph on $n$ vertices with parts of size $\left\lceil \frac{n}{t}\right\rceil$ or $\left\lfloor \frac{n}{t}\right\rfloor$) has the most edges among all $K_{t+1}$-free graphs on $n$ vertices.
Erd\H os \cite{MR0329905} and Bollob\'as, Erd\H os, and Szemer\'edi~\cite{BES75} asked the following Tur\'an-type problem for multipartite graphs.
\begin{problem}
\label{pro:BES}
Given integers $n$ and $2\le t\le r-1$, what is the largest minimum degree $\delta(G)$ among all $r$-partite graphs $G$ with parts of size $n$ and which do not contain a copy of $K_{t+1}$?
\end{problem}

Let $f(n, r, t+1)$ denote the answer to Problem~\ref{pro:BES}. At a meeting in 1972, Erd\H os conjectured that $f(n, r, r)= (r-2)n$, see \cite[Problem 2, 353--354]{MR0329905}. Graver gave a short and elegant proof for $r=3$ but Seymour constructed counterexamples for $r\ge 4$, see~\cite{BES75}. The study of $f(n, r, r)$ (mostly in its complementary form concerning \emph{independent transversals})  has been a central topic in Combinatorics (see, e.g.,~\cite{MR955135, MR2195582, MR1860440, MR1208805, MR2246152}) due to its applications in graph arboricity, list coloring, and strong chromatic numbers. The problem of determining $f(n, r, r)$ was finally settled by Haxell and Szab\'{o} \cite{MR2195582}, and Szab\'{o} and Tardos \cite{MR2246152}; indeed, for every $n\in \mathbb{N}$ and even $r\ge 2$,
\begin{align} \label{eq:rr}
f(n, r+1, r+1)-n= f(n, r, r) = (r-1)n - \left\lceil \frac{rn}{2(r-1)}\right\rceil.
\end{align}

In contrast, little is known about the value of $f(n, r, t+1)$ for $r> t+1$. In 1975 Bollob\'as, Erd\H os, and Szemer\'edi~\cite{BES75} stated Problem~\ref{pro:BES} explicitly
and noted that Tur\'an's theorem easily implies that 
\begin{align}
    \label{eq:BES}
    f(n, r, t+1)= \left(r - \frac{r}{t}\right) n \quad \text{when $t$ divides $r$}. 
\end{align}
Indeed, for any $r\ge t+1$, Turan's theorem implies that every $K_{t+1}$-free graph $G$ on $rn$ vertices has at most $(1- 1/t)(rn)^2/2$ edges, and thus $\d(G)\le (1- 1/t)rn$. 
On the other hand, we may let $G$ be the complete $t$-partite graph on $r n$ vertices with parts of size $\left\lceil \frac{r}{t}\right\rceil n$ or $\left\lfloor \frac{r}{t}\right\rfloor n$ (in other words, $G$ is an $n$-vertex blow-up of the Tur\'an graph $T_{t}(r)$). Then 
\begin{align}\label{eq:TuranBd}
\left( r - \left\lceil \frac{r}{t}\right\rceil \right) n \le f(n, r, t+1)\le  \left( r -  \frac{r}{t} \right) n.
\end{align}
%We will also study the $r>t$ case of Problem~\ref{pro:BES} in its complement form. 
%Let $\Delta(n, r, t) := (r-1)n - f(n, r, t)$ denote the smallest $\Delta(G)$ among all  $r$-partite graphs $G$ with parts $V_1,\dots, V_r$  of size $n$ and without a crossing independent  of size $t$ (i.e., an independent set of size $t$ with at most one vertex in each $V_i$).
%While studying the Tur\'an number of disjoint copies of $K_t$ in $r$-partite graphs, Han and Zhao~\cite{HZ19} encountered the problem of $f(n, t, r)$ on not necessarily balanced $r$-partite graphs. 
%Unaware of known constructions from~\cite{BES75, MR2246152}, they made a conjecture~\cite[Conjecture 4.1]{HZ19} which generalizes \eqref{eq:TuranBd} for not necessarily balanced $r$-partite graphs.\COMMENT{AT: is this needed now as the conjecture looks like it is deleted in the paper}
%We know that this conjecture fails when $r=t\ge 4$ but it was shown to be true when $t=3$.
Extending Graver's work on $f(n, 3, 3)$, Bollob\'as, Erd\H{o}s, and Straus~\cite{BEST} answered Problem~\ref{pro:BES} for all (not necessarily balanced) $r$-partite graphs $G$ when $t=2$. Their result implies that for every $n\in \mathbb{N}$ and $r \geq 3$, 
\[
f(n,r,3)= \lfloor r/2 \rfloor n.
\]
 
% \begin{theorem} 
% \label{thm:HZ}\cite{BEST}
% Let $G$ be an $r$-partite graph whose parts have sizes $n_1,\dots, n_r$. Let $I, J$ be an optimal bipartition of $[r]$, namely, a bipartition that minimizes $| n_I - n_J |$ among all bipartitions of $[r]$ (where $n_I := \sum_{i\in I} n_i$). If $G$ contains no triangle, then $\delta(G)\le \min\{ n_I, n_J \}$. 
% \end{theorem}

\medskip
The aim of this paper is to rebuild momentum on Problem~\ref{pro:BES} for $r> t+1\geq 4$.
For any such choice of $r$ and $t$, our results either resolve Problem~\ref{pro:BES} or provide a lower bound on $f(n,r,t+1)$ that improves that given in 
 (\ref{eq:TuranBd}). 
 In particular, 
 in the case that $r \equiv -1 \pmod{t}$, our first result shows that 
  the lower bound in \eqref{eq:TuranBd} is tight.

\begin{theorem}\label{thm:-1}
Given integers $n\ge 1$, $m\ge 2$, and $t \ge 3$, let $r=mt-1$. Then
$f(n, r, t+1) =\left( r - \left\lceil {r}/{t}\right\rceil \right) n$. %= (mt-1-m)n
\end{theorem}
It turns out that the lower bound in \eqref{eq:TuranBd} is best possible only if
$r \equiv 0,-1 \pmod{t}$; in all other cases we give constructions that improve on this lower bound (see Section~\ref{sec:construct}). Moreover, 
in many such cases, including when $r \ge (3t-1)(t-1)$, we determine $f(n,r,t+1)$ up to an additive constant.
\begin{theorem}~\label{thm:-2}
Given integers $n\ge 1, m\ge 2, t \ge 3$, let $r = mt -a$ with $2 \le a \le \min\{m, t-1\}$. Suppose
%of the following conditions holds,
\begin{itemize}
    \item[(i)] $r \ge a(3t-1)$  or 
    \item[(ii)] $\frac{r}{t(3t-1)(m-1)}-\frac{a}{t(m-1)} + \frac{a-1}{mt-2} \geq  \frac{1}{n} $.
\end{itemize}
Then
\begin{align*}
(r-1)n - (m-1) \left\lceil \frac{(r-1)n}{mt-2} \right \rceil \le f(n,r,t+1) \leq 
(r-1)n - \left\lceil\frac{(m-1)(r-1)n}{mt-2}\right\rceil
.
\end{align*}
% Suppose that
% \begin{align*}
% \frac{r}{t(3t-1)}-\frac{a}t + \frac{(m-1)(a-1)}{mt-2} >  0 
% \end{align*} 
% and $n$ is sufficiently large. \COMMENT{AT: need the sufficiently large here as we apply Cor~\ref{newcor}.
% Or we could instead apply Cor~\ref{cor:AES1} and get a better bound on $n$ at the expense of not having $\frac{r}{t(3t-1)}-\frac{a}t + \frac{(m-1)(a-1)}{tm-2} >  0 $ in the statement of the theorem.  
% AL: added a moreover statement to mention to latter thing.}
% Then
% \al{Added the moreover statement, which does not need $n$ sufficiently large.}Moreover, if $r \ge a(3t-1)$ then \eqref{eqn:main} holds for all $n \in \mathbb{N}$. 
\end{theorem}
%It is easy to see that the full hypothesis of Theorem~\ref{thm:-2} holds when $2\leq a \leq t-1$ and $r\geq (3t-1)a$. 
Thus, up to an additive constant, \eqref{eq:BES}, Theorems~\ref{thm:-1} and~\ref{thm:-2} together resolve Problem~\ref{pro:BES} whenever
%for all choices of $r,t$ satisfying 
$r \geq (3t-1)(t-1)$. In particular, 
%if $r= mt -a$ with $a\ge 1$, $m\ge 3a$, and 
if $r \geq (3t-1)(t-1)$, $r \not\equiv 0 \pmod{t}$, and 
$\lceil r/t \rceil t-2$ divides $n$, then
\[
f(n, r, t+1) = (r-1)n - \frac{\lceil r/t \rceil-1}{\lceil r/t \rceil t-2} (r-1)n.
\]

Combining the last two theorems with~\eqref{eq:BES}, we essentially determine $f(n,r,4)$ for all $r \not = 7$. 
\begin{corollary}~\label{thm:t4}
Let $r \ge 5$ with $r \ne 7$. Suppose $n \geq 60$ if $r=10$; $ n \geq 22$ if $r=13$, and $n\in \mathbb N$ otherwise.\COMMENT{YZ rephrased this sentence.} Then 
\begin{align*}
	f(n,r,4) = 
	\begin{cases}
	\left(r- \lfloor r/3 \rfloor \right)\frac{r-1}{r} n - c_r & \text{if $r \equiv 1 \pmod{3}$,}\\
	\left( r - \left\lceil {r}/{3}\right\rceil \right) n &\text{otherwise},
	\end{cases}
\end{align*}
where $0 \le  c_r \le r/3 $.
\end{corollary}
\al{For $f(n,7,4)$, I know that it is between $30n/7= 4.28 n $ and $4.31 n $. I can deduce that there is an independent of size at least $2.46 n$ which spans exactly 3 vertex classes. }

The proof of Theorem~\ref{thm:-2} applies a result of Andr\'asfai, Erd{\H{o}}s, and S\'os~\cite{MR0340075}. In particular, this result allows us conclude that, if $r$ is large compared to $t$, then 
 $f(n,r,t+1)$ is equal to the maximum of the minimum degrees among the family of balanced $r$-partite $rn$-vertex graphs of chromatic number at most $t$.
 This approach is also utilised in the proof of Theorem~\ref{thm:-1} when $m \geq 3$. Interestingly, this approach breaks down when $m =2$, so we  require a separate direct argument in this case (which hinges on the fact that $r=2t-1$).

\smallskip

There are two extremal problems closely related to Problem~\ref{pro:BES}. First, a multipartite Tur\'an theorem has been known since the 1970s. Bollob\'as, Erd\H os, and Szemer\'edi~\cite{BES75} showed that the $n$-vertex blow-up of $T_t(r)$ has the most edges among all $r$-partite $K_{t+1}$-free  graphs with $n$ vertices in each part. In fact, this easily follows from Tur\'an's theorem and indicates that for multipartite graphs, the minimum degree version, Problem~\ref{pro:BES}, is indeed harder than the Tur\'an problem.\footnote{In contrast, Tur\'an's theorem easily implies that $\max \d(G)= n - \lceil n/t \rceil$ among all $K_{t+1}$-free graphs on $n$ vertices.} 
Furthermore, Bollob\'as, Erd\H{o}s, and Straus~\cite{BEST} determined the largest size of (not necessarily balanced) $r$-partite $K_{t+1}$-free graphs. Another related problem concerns finding the smallest $d_{r}^{t}$ such that every $r$-partite graph whose parts have pairwise edge density greater than $d_{r}^{ t}$ contains a copy of $K_{t}$. Bondy, Shen, Thomass\'e, and Thomassen~\cite{bondy} showed that $d_{3}^{3}=0.618...$, the golden ratio. Pfender~\cite{pfender} showed that $d_r^t= (t-2)/(t-1)$ for sufficiently large $r$; in particular, $d_r^3= 1/2$ for $r\ge 13$.

\subsection{Notation}
Given a graph $G$ and $x \in V(G)$ we write $N(x)$ for \emph{the neighbourhood of $x$ in $G$} and define $d (x):=|N(x)|$ as the \emph{degree of $x$} in $G$. If $X \subseteq V(G)$ we write $N(x,X):=N(x)\cap X$ and 
$d (x,X) :=|N(x,X)|$. We write $G[X]$ for the induced subgraph of $G$ with vertex set $X$.

Let $G$ be an $r$-partite graph with parts $V_1,\dots, V_r$.
A \emph{transversal} is a subset $S\subset V(G)$ so that $|S\cap V_i|=1$ for each part $V_i$ of $G$.
 A set $S\subset V(G)$ is \emph{crossing} if $|S\cap V_i|\le 1$ for each part $V_i$ of $G$. 
 %We say a subgraph $H$ of $G$ is \emph{crossing} if $V(H)$ is crossing.
 \COMMENT{YZ: why do we need this sentence? And why don't we use independent transverals for the rest of the paper? {\bf AT: I introduced this as I found it confusing that the words `independent transveral' means a `spanning' independent set but an `independent transversal of size X' means something a bit different. So I'd rather keep this in to avoid any reader confusion in Section 4. In particular, at one point during section 4 we need to say `crossing independent set' without specifying its size. If we replace it by just `independent transversal' the reader would think we meant spanning when we really don't. YZ: agree.}\\
 Also, the  sentence `We say a subgraph $H$ of $G$ is \emph{crossing} if $V(H)$ is crossing.' was added as we talk about a crossing independent set $\bar{K_t}$ so formally this is a subgraph rather than subset... but don't mind if we delete this. \textbf{YZ has removed this sentence and revised all the sentences that we had used $\overline{K}_t$}} 
 Thus an independent transversal is simply a crossing independent set of size $r$.
Define $K_r(n)$ to be the complete $r$-partite graph where each part  has size~$n$.

\subsection{Organization of paper}
\al{\textbf{Edited this section due to organising the sections.}}
In Section~\ref{sec:comp} we formally restate Problem~\ref{pro:BES} in its complementary form and prove Theorem~\ref{thm:-1} for $m=2$.
In Section~\ref{sec:AES}, we introduce a related parameter $\delta(n,r,t)$ that is equal to  $f(n,r,t+1)$ if $r$ is large compared to~$t$.
 In Section~\ref{sec:construct} we give constructions that improve on the lower bound in \eqref{eq:TuranBd} whenever
$r \not \equiv 0,-1 \pmod{t}$.
%In Section~\ref{sec:f(n,2t-1,t+1)}, we prove Theorem~\ref{thm:-1} for $m=2$. 
%evaluate $f(n,2t-1,t+1)$ for $t \geq 3$.
We give an upper bound on $\delta(n,r,t)$ in Section~\ref{sec:bounds} that allows us to easily deduce 
Theorems~\ref{thm:-1} and~\ref{thm:-2} in Section~\ref{sec:main}.  We finish the paper with concluding remarks in Section~\ref{sec:conclude}.

\section{The complementary problem and Proof of Theorem~\ref{thm:-1} for $m=2$}% and a related parameter}
%\subsection{The complementary problem}
\label{sec:comp}
In most papers on $f(n, r, r)$, the complementary form of $f(n, r, r)$ was considered, namely, the smallest $\Delta(G)$ among all $r$-partite graphs $G$ with parts $V_1,\dots , V_r$ of size $n$ and without an independent transversal. For a couple of our proofs it will also be easier to work in the corresponding complementary setting also.

 Let $\Delta(n, r, t) := (r-1)n - f(n, r, t)$ denote the smallest maximum degree $\Delta(G)$ among all  $r$-partite graphs $G$ with parts $V_1,\dots, V_r$  of size $n$ and without a crossing independent set  of size $t$.
Note that
\eqref{eq:TuranBd} is equivalent to \al{\textbf{removed labelling in the mathdisplay}}
\begin{align*}
\left( \frac{r}{t} - 1 \right) n \le \Delta(n, r, t+1) \le  \left( \left\lceil \frac{r}{t}\right\rceil - 1 \right) n.
\end{align*}

%
%\subsection{Proof of Theorem~\ref{thm:-1} for $m=2$} \label{sec:f(n,2t-1,t+1)}

%In this section we prove the following lemma,
%give a direct proof that $f(n,2t-1,t+1) = (2t-3)n$ for $t\ge 3$. This 
\al{\bf modified linking sentence.} We now prove the following lemma, 
which is the $m=2$ case of Theorem~\ref{thm:-1}. 
\begin{lemma} \label{lma:f(n,2t-1,t+1)}
For every $n\in \mathbb{N}$ and $t \ge 3$, $f(n,2t-1,t+1) = (2t-3)n$.
\end{lemma}
\begin{proof}
Let $r := 2t-1$.
By~\eqref{eq:TuranBd}, it suffices to prove that $f(n,2t-1,t+1) \le (2t-3)n$. Further,
by considering the complementary problem, it suffices to prove that if $G$ is an $r$-partite graph with vertex classes $V_1, \dots, V_r$ of size $n$ so that $G$ does not contain a crossing independent set of size $t+1$, then $\Delta(G) \ge n $. 

Suppose for a contradiction that $\Delta(G) < n$.

\begin{claim} \label{clm:new}
For all $x \in V(G)$ and $i \in [r]$, we have $d(x,V_i) < n/(t-1)$.
\end{claim}

\begin{proofclaim}
Suppose the claim is false.
Let $D:=\max\{ d(v,V_i) \colon v \in V(G), i \in [r]\}\geq n/(t-1)$.
Without loss of generality, we may assume that there exists $x_1 \in V_1$ such that $d(x_1,V_2) =D$.
Since $\Delta (G) <n$, there exists $x_2 \in V_2 \setminus N(x_1)$. Furthermore, since $\Delta(G) < n$ and $r =2t-1$, we can greedily find $x_3, \dots, x_{t-1}$ such that $S = \{x_1,\dots, x_{t-1}\}$ is a crossing independent set.\footnote{We can actually find a crossing independent set of size $t$. However, considering a crossing independent set of size $t-1$ is crucial to the argument here.}
Without loss of generality, assume that $x_i \in V_i$ for $i \in [t-1]$.

For each $i\in [ t]$, let $W_i$ consist of all the vertices of $V_{r+1-i}$ that are not adjacent to any vertex in~$S$.
Set $n_i := |W_i|$ and without loss of generality, assume that $n_1 \ge n_2 \ge \dots \ge n_{t}$.
Let $\ell := \max \{ i \in [t] \colon n_i > 0 \}$. 
Note that
\begin{align*}
    n_1+ \dots + n_{\ell} 
    & = n_1 + \dots +n_{t} 
    = \left| \bigcup_{t \le i \le r} V_{i} \setminus \bigcup_{j \in [t-1]}N(x_j) \right| \\
    & \ge tn - ((t-1)\Delta(G)-d(x_1,V_2)) > n+D
.
\end{align*}
By averaging, we have 
\begin{align}  \label{eqn:n_1+n_l-1}
        n_1+ \dots + n_{\ell-1} > (\ell-1) (n+D) / \ell. 
\end{align}
Notice that the $\ell$-partite subgraph of $G$ induced by $W_1,\dots, W_{\ell}$ must be complete (otherwise one can extend~$S$ into a crossing independent  set of size~$t+1$, a contradiction).
Hence, there exists $y \in W_{\ell}$ such that $\bigcup_{i \in [\ell-1]} W_i \subseteq N(y)$.
By the definition of $D$, we deduce that 
\begin{align*}
    n_1 + \dots + n_{\ell-1} \le \sum_{i \in [\ell-1]} d(y,V_{r+1-i})  \le (\ell-1)D.
\end{align*}
Together with~\eqref{eqn:n_1+n_l-1}, this implies that $D > n/(\ell-1)$ and so 
\begin{align*}
    \Delta(G) \ge d(y) \ge n_1 + \dots + n_{\ell-1} >  (\ell-1) (n+D) / \ell \ge n,
\end{align*} 
a contradiction. 
\end{proofclaim}

Given a crossing independent set~$S$ of~$G$, let $\sigma(S) := \sum_{x \in S} d(x,V_S)$ where $V_S$ in the union of all $V_i$ that contain a vertex from $S$.
As mentioned earlier in the proof, we can greedily construct a crossing independent set of size~$t$.
Let $S$ be a crossing independent  set of size~$t$ with $\sigma (S)$ maximal. 
Without loss of generality, $S \cap \bigcup_{i \in [t-1]} V_i= \emptyset$.

Consider any $(t-1)$-set $S'\subset S$.
By Claim~\ref{clm:new}, there exists a vertex $y \in V_1$ that is not adjacent to any vertex in~$S'$.
Note that $S'\cup \{y\}$ is a crossing independent set of size~$t$.
Hence $	\sigma(S)  \ge  \sigma(S'\cup \{y\} ) \ge \sigma(S') + \sum _{x\in S'} d(x,V_1)$.
By summing over all  $(t-1)$-sets $S'\subset S$, we obtain that
\begin{align*}
	t \sigma(S) & \ge \sum_{S' \subset S \colon |S'|=t-1} \sigma(S') + (t-1) \sum_{x \in S}d(x,V_1) 
	= (t-2)\sigma(S)  + (t-1) \sum_{x \in S}d(x,V_1)\\
	& \ge (t-2)\sigma(S) + (t-1)n,
\end{align*}
where the last inequality follows as every vertex in $V_1$ must be adjacent to at least one vertex in $S$ (else there exists a crossing independent  set of size~$t+1$, a contradiction).
Thus $ \sigma(S) \ge (t-1) n/2$.

As $S \cap \bigcup_{i \in [t-1]} V_i= \emptyset$, then $\bigcup_{i \in [t-1]} V_i \subseteq \bigcup_{x \in S} N(x)$ or else there exists a crossing independent set of size~$t+1$, a contradiction. 
Therefore,
\begin{align*}
t \Delta(G) 
 \ge \sum_{x \in S} d(x) \ge \sigma(S) + \left|\bigcup_{i \in [t-1]} V_i \right|
 \ge \frac{(t-1)n}2 + (t-1) n  \ge t n,
\end{align*}
implying that $\Delta(G) \ge n$, a contradiction.
\end{proof}

%%%%%%%%%%%%%%%%%%%%%%%%%%%%%%%%%%%%%%%%%%%%%%%%%%%%%%%%%%%%%%%%%%%%%%%%%%%%%

\section{A connection to the parameter $\delta(n,r,t)$ } \label{sec:AES}

The following problem turns out to be closely related to Problem~\ref{pro:BES}.
Let $\mathcal{G}(n,r,t)$ be the family of all $r$-partite graphs $G$ with parts of size~$n$ and with chromatic number $\chi(G) \le t$.
Let $\delta(n,r,t) := \max \{ \delta(G) \colon G \in \mathcal{G}(n,r,t)\}$.
An $n$-vertex blow-up of the Tur\'an graph $T_{t}(r)$ is a member of $\mathcal{G}(n,r,t)$. Together with \eqref{eq:TuranBd}, this gives 
\begin{align} \label{eq:TuranBd2}
	\left( r - \left\lceil \frac{r}{t}\right\rceil \right) n \le \delta(n,r,t) \le f(n, r, t+1)\le  \left( r -  \frac{r}{t} \right) n.
\end{align}
%On the other hand, we have \[\d(n, r, t)\le \left(1 - \frac{1}{t}\right) r n \] because every graph on $rn$ vertices with $\chi(G)\le t$ satisfies $\alpha(G)\ge rn/t$ and $\d(G)\le rn - rn/t$. 
Therefore, when $t$ divides $r$, we have $\d(n,r,t) = f(n, r, t+1) = (r - r/t)n$.

When $r$ is large compared to $t$, Corollary~\ref{cor:AES1} below implies that $f(n,r,t+1) = \delta(n,r,t)$ as well. In fact, this is an easy consequence of the following result of Andr\'asfai, Erd\H{o}s, and S\'os~\cite{MR0340075}. \COMMENT{AT: observe that it isn't always the case that extremal e.g.s are $(t-1)$-partite... e.g. $4$-partite case for $K_4$... double check. YZ, I don't understand the question. Are you referring to Theorem 2.1 or Corollary 2.3? {\bf AT: what I am really asking here is do we know of cases when $\d(n,r,t) \ne f(n, r, t+1)$?}}

\begin{theorem}[Andr\'asfai, Erd{\H{o}}s, and S\'os~\cite{MR0340075}] \label{thm:AES}
Let $t \ge 2$ and let $G$ be a $K_{t+1}$-free graph on $N$ vertices.
If $\delta(G) > \frac{3t-4}{3t-1} N$, then $\chi (G) \le t$. 
\end{theorem}

%We obtain the following immediate corollaries. 

\begin{corollary}\label{cor:AES}
For $r > t \ge 2$,
%If $\delta(n,r,t) \ge \frac{3t-4}{3t-1} rn$, then $f(n,r,t+1) = \delta(n,r,t)$.
$f(n,r,t+1) \le \max \left\{ \frac{3t-4}{3t-1} (rn), \delta(n,r,t) \right\}$.
\end{corollary}
\begin{proof}
Let $G$ be a $K_{t+1}$-free $r$-partite graph with $n$ vertices in each part. If $\delta(G)> \frac{3t-4}{3t-1} (rn)$, then $\chi(G)\le t$ by Theorem~\ref{thm:AES}. Thus $G\in \mathcal{G}(n, r, t)$ and $\delta(G)\le \delta(n, r, t)$.
\end{proof}
By applying~\eqref{eq:TuranBd2} together with Corollary~\ref{cor:AES} one can 
conclude that $f(n,r,t+1) = \delta(n,r,t)$ provided that $r\geq (t-1)(3t-1)$.
 
\begin{corollary}\label{cor:AES1}
Let $r > t \ge 2$ and $0\leq a \leq t-1$ so that $r \equiv -a \pmod{t}$. If $ r \ge a(3t-1)$
then $f(n,r,t+1) = \delta(n,r,t)$.
\end{corollary}

\begin{proof} 
\eqref{eq:BES} covers the case when $a =0$ so we assume that $a \in [t-1]$.
By~Corollary~\ref{cor:AES}, it suffices to show that $\delta(n,r,t) \ge \frac{3t-4}{3t-1} rn$.
By~\eqref{eq:TuranBd2}, we have 
\begin{align*}
\delta(n,r,t) &\ge \left( r - \left\lceil \frac{r}{t}\right\rceil \right) n = \left( r - \frac{r+a}{t} \right) n 
= \left( \frac{r}{t(3t-1)} + \frac{3t-4}{3t-1}r -\frac{a}{t} \right) n \ge  \frac{3t-4}{3t-1} rn,
\end{align*}
where we use the fact that $r \ge a(3t-1)$ in the last inequality.
\end{proof}

%%%%%%
\section{Lower bound constructions}\label{sec:construct}
 In this section we give constructions that improve on the lower bound in \eqref{eq:TuranBd} whenever
$r \not \equiv 0,-1 \pmod{t}$.
\begin{proposition} \label{prop:d(n,r,t)3lower}
Let $r \ge m,t \ge 2$ be such that $m (t-1) \le r \le mt-1$.
Then there exists a graph $G \in \mathcal{G}(n,r,t)$ such that $\delta(G) = (r-1)n - (m-1) \left\lceil \frac{(r-1)n}{mt-2} \right\rceil$.
In particular,
\begin{align*}
   f(n,r,t+1) \geq  \delta(n,r,t) \ge 	(r-1)n - (m-1) \left\lceil \frac{(r-1)n}{mt-2} \right\rceil.
\end{align*}
%In other words, $\Delta(n,r,t+1) \le (m-1) \left\lceil \frac{(r-1)n}{mt-2} \right\rceil$.
\end{proposition}

\begin{proof}
When $r=mt-1$, the desired bound follows from \eqref{eq:TuranBd2}. We thus assume $r\le mt-2$.
Let $\ell := \lceil (r-1)n/(mt-2) \rceil \le n $.
Let $K:=K_r(n)$ and let  $V_1, \dots , V_{r}$ denote its parts. 
For $i \in [t-1]$, let $B_i := \{(i-1)m+1, \ldots, im\}$.
So $B_1, \dots, B_{t-1}$ form an equipartition of~$[m(t-1)]$.
For  $i \in [t-1]$, let $W_i \subset \bigcup_{j \in B_i} V_j$ be such that $|W_i \cap V_j| = \ell$ for $j \in B_i$.
Let $W_t := V(K) \setminus \bigcup_{i \in [t-1]}W_i$. 
Then 
\[
|W_1| = \cdots = |W_{t-1}|= \ell m \quad \text{and} \quad |W_t|= rn - m (t-1) \ell .
\]
Let $G'$ be the complete $t$-partite graph with parts $W_1, \dots , W_{t}$.
We set $G := K \cap G'$; that is, $G$ is the graph on $V(K)$ such that, for $x \in V_i \cap W_j$ and $x' \in V_{i'} \cap W_{j'}$, we have $xx' \in E(G)$ if and only if $i \ne i'$ and $j \ne j'$.
Clearly $\chi(G) \le \chi (G') \le t$. 
If $x \in V_i$ with $i \notin [m(t-1)]$, then $d(x) = r n - |W_t|$.
If $x \in V_i$ with $i \in [m(t-1)]$, then 
\begin{align*}
    d(x) = 
    \begin{cases}
     rn -  (n + |W_1| - \ell )  
     %= (r-1)n - (m-1)\lceil \alpha n  \rceil
     & \text{if $x \notin W_t$,}\\
     rn -  (n + |W_t| - ( n- \ell ) )  
     %=  (r-1) n - \left( (m-2)+ (m(t-1)-1) \lfloor (1-\alpha)n \rfloor \right)
     & \text{if $x \in W_t$.}    \end{cases}
\end{align*}
By our choice of $\ell$,
\begin{align*}
	\delta(G) &  = (r-1) n -  \max \{ (m-1) \ell,  (r-1)n - (m(t-1)-1) \ell  \}\\
	& =  (r-1) n - (m-1) \ell 
	=  (r-1)n - (m-1) \left\lceil \frac{(r-1)n}{mt-2} \right\rceil ,
\end{align*}
as required.
\end{proof}
Note that the lower bound in Proposition~\ref{prop:d(n,r,t)3lower} improves the lower bound from (\ref{eq:TuranBd}) in the case when $r=mt-a$ with $2\leq a \leq \min\{m, t-1\}$.
Indeed, in this case (\ref{eq:TuranBd}) gives a lower bound of $(r-m)n$ while 
\begin{align}
(r-1)n - (m-1) \left\lceil \frac{(r-1)n}{mt-2} \right\rceil &> (r-1)n - (m-1) \frac{(r-1)}{mt-2}n-(m-1) \nonumber \\
&= (r-m)n + \frac{(a-1)(m-1)}{mt-2}n - (m-1).
\label{eq:LB}
\end{align}
Thus, if $n\ge (mt-2)/(a-1)$, then the lower bound in Proposition~\ref{prop:d(n,r,t)3lower} improves the lower bound from (\ref{eq:TuranBd}).

In the remaining case -- when $m<a\leq t-1$ -- the next result beats the lower bound from (\ref{eq:TuranBd}) when $n$ is not too small.

\begin{proposition} \label{prop:construction}
Let $r > t \geq 3$ be such that $r = mt - a$ with $2\leq m <a<t $.
Then 
\begin{align*}
	 f(n,r,t+1) \geq \delta(n, r, t) \ge  (r-1)n - (m-1) \left\lceil \frac{(m(t-1-a+m)-1)n}{m(t-a+m)-2} \right\rceil.
\end{align*}

\end{proposition}

\begin{proof}
Let $t': = t-a+m$ and $r' := m(t'-1)$. 
By Proposition~\ref{prop:d(n,r,t)3lower}, there exists a graph  $G' \in\mathcal{G}(n,r',t')$ such that 
\begin{align*}
\delta(G') = (r'-1)n -  (m-1) \left\lceil \frac{(r'-1)n}{mt'-2} \right\rceil 
\le (r'-1)n - (m-2)n,
\end{align*}
where the last inequality is due to fact that 
\begin{align*}
    (m-1)(r'-1) - (mt'-2)(m-2) = m(t'-m+2)-3  = (t-a+2)m-3 > 0. 
\end{align*}

We now construct a graph $G \in \mathcal{G}(n,r,t)$ from $G'$ as follows. 
Let $W_{a-m+1} := V(G')$.
Let $W_1, \dots, W_{a-m}$ be vertex sets each of size $(m-1)n$ such that $W_1, \dots, W_{a-m+1}$ are disjoint. 
Let $G$ be the resulting graph on $\bigcup_{i \in [a-m+1]} W_i$ obtained from $G'$ by adding edges $xx'$ for all $x \in W_i$ and $x' \in W_{i'}$ with $i \ne i'$.
Note that $\chi(G) = \chi(G') + a-m \le t$. 
Since each $W_i$ with $i \in [a-m]$ can be partitioned into $m-1$ vertex classes of size $n$, we deduce that $G \in \mathcal{G}(n,r,t)$. 

For $x \in V(G) \setminus V(G')$, we have  $d_G(x) = r n -(m-1) n $,  and for $x \in V(G')$ we have
\begin{align*}
    d_G(x) = (a-m)(m-1) n + \delta(G') =  (r-1)n -  (m-1) \left\lceil \frac{(r'-1)n}{mt'-2} \right\rceil 
    \le  (r-1)n - (m-2)n.
\end{align*}
Hence,
\begin{align*}
    \delta(G)  =  (r-1)n -  (m-1) \left\lceil \frac{(r'-1)n}{mt'-2} \right\rceil,
\end{align*}
as required. 
\end{proof}

%\al{this is better than the other bound of $-(m-1)(r-1)/(r-m+2)$ for $r = (m-1)t -b$. I take $(a = m+b)$, so
%\begin{align*}
%     \frac{(r'-1)n}{mt'-2} - \frac{r-1}{r-m+2} = %\frac{-b(m^2-3m+1)-(mt-2)}{(r-m+2)(mt'-2)}<0
%\end{align*}
%Please check!
%}
%Actually, if we combine with Proposition~\ref{prop:d(n,r,t)3} (see later), our construction actually shows that $ \delta(n,r',t') \ge (a-m)(m-1) n + \delta(n,r',t')$ for large enough~$n$. 

By applying Corollary~\ref{cor:AES} with Proposition~\ref{prop:d(n,r,t)3lower}  we can obtain the following result, which improves on Corollary~\ref{cor:AES1} in most cases when $n$ is not too small.

\begin{corollary}\label{newcor}
Given $m, t \ge 2$, let $r = mt -a$ where $2 \le a \le \min\{m, t-1\}$.\COMMENT{\bf YZ added $a\le t-1$ because the longest displayed line in the proof requires that $m-\lceil r/t\rceil$.}
If
\begin{align*} %\label{gbound}
\frac{r}{t(3t-1)(m-1)}-\frac{a}{t(m-1)} + \frac{a-1}{tm-2} \geq  \frac{1}{n},
\end{align*} 
\COMMENT{{\bf AT: I think this might be the best way to say it. That is, rather than decoupling the conditions, this condition allows $n$ to be smaller as you make $r$ bigger}. YZ: agree.}
then $f(n,r,t+1) = \delta(n,r,t)$.
\end{corollary}
%Given any $r \ge t \ge 2$ and $1\leq a \leq t-1$\COMMENT{AT: $a \geq 2$ or $r>a(3t-1)$} such that
%$r \equiv -a \pmod{t}$ and $ r \ge a(3t-1)$, notice that $r$ can be written in the form
%$r = mt -a$ where  $m \geq 2$ and $a \le m/3$, and so that (\ref{gbound}) holds. Thus $f(n,r,t+1) = \delta(n,r,t)$ in this case.
%Note that (\ref{gbound}) holds when $r \geq a(3t-1)$ provided $n \geq \frac{tm-2}{a-1}$.

\begin{proof}
By (\ref{eq:TuranBd2}) and Corollary~\ref{cor:AES} it suffices to prove that 
\begin{align}\label{aim}
    \delta(n,r,t) \geq \frac{3t-4}{3t-1} rn.
\end{align}
Proposition~\ref{prop:d(n,r,t)3lower} and \eqref{eq:LB} together imply that
\[
\delta(n,r,t)> (r-m)n + \frac{(a-1)(m-1)}{mt-2}n - (m-1). 
\]
On the other hand, $\frac{3t-4}{3t-1} r = r-m +\frac{a}t-\frac{r}{t(3t-1)}.$
Thus, to prove \eqref{aim}, it suffices to have
\begin{align*}
\frac{(a-1)(m-1)}{mt-2}n - (m-1)\ge 
%+
\frac{an}t-\frac{rn}{t(3t-1)}.
%    (r-1)n - (m-1) \frac{(r-1)}{mt-2}n-(m-1) \geq \frac{3t-4}{3t-1} rn = rn-mn +\frac{an}t-\frac{rn}{t(3t-1)}.
\end{align*}
%As in the proof of Corollary~\ref{cor:AES1},
%$$r-m=r-\left \lceil \frac{r}{t} \right \rceil=\frac{r}{t(3t-1)}-\frac{a}{t}+ \frac{3t-4}{3t-1} r.$$
Indeed, this is equivalent to our assumption
\[
\frac{r}{t(3t-1)(m-1)}-\frac{a}{t(m-1)} + \frac{a-1}{mt-2} \geq  \frac{1}{n}. \qedhere
\]
\end{proof}

%%%%%%%%%%%%%%%%%%%%%%%%%%%%%%%%%%%%%%%%%%%%%%%%%%%%%%%%%%%%%%%%%%%%%%%%%%%%

\section{An upper bound on $\delta(n,r,t)$} \label{sec:bounds}

In this section we prove the following upper bound on $	\delta(n,r,t)$.
%%%%%
\begin{proposition} \label{prop:d(n,r,t)3}
Let $r,n \in \mathbb{N}$ and $m,t \ge 2$ be such that $(m-1)t < r < mt$.
Then,
\begin{align*}
	\delta(n,r,t) \le (r-1)n - \left\lceil \frac{(m-1)(r-1)n}{mt-2}\right\rceil.
\end{align*}
%In particular, for $t=3$, we have  $\delta(n,3p+1,3) \ge 3p(2p+1) n /(3p+1) - p $.
\end{proposition}

\begin{proof}
Let $\Delta^* :=  (m-1) (r-1)n/(mt-2) $. Since $\delta(n,r,t)$ is an integer, it suffices to show that $\delta(n,r,t) \le (r-1)n - \Delta^*$.
Suppose to the contrary that $\delta(n,r, t) > (r-1)n -  \Delta^* $.
Let $G \in \mathcal{G}(n,r,t)$ be such that $\delta(G) =\delta(n,r, t)$. 
As $G \in \mathcal{G}(n,r,t)$, $V(G)$ can be partitioned into $r$ independent sets $V_1,\dots , V_r$ each of size $n$; 
as $\chi(G) \le t$, $V(G)$ can be partitioned into $t$ color classes $W_1, \dots, W_{t}$.
\COMMENT{YZ removed ``After adding edges if necessary, we may assume that, whenever $x\in V_i\cap W_j$ and $x'\in V_{i'}\cap W_{j'}$, $xx'\in E(G)$ if and only $i\ne i'$ and $j\ne j'$.
In other words, all the edges between two color classes in two different parts are present''. I remember the purpose of having them was for a claim that we no longer need.
{\bf AT: replaced some equalities with inequalities in proof to make it work}} %We call a 2-dimensional array $\mathbf{w}:= \{w_{i,j}\}_{i\in [r], j\in [t]}$ \emph{feasible} if all $w_{i, j}$ are non-negative integers and $\sum_{j=1}^t w_{i, j} = n$ for all $i\in [r]$.
%Thus $G = G_{\mathbf{w}}$ is uniquely determined by some feasible $\mathbf{w}$.
%In particular, $e(G)= \sum_{i\ne i', j\ne j'} w_{i, j} w_{i', j'}$.
For every $x\in V_i\cap W_j$, we have $d(x)\leq(r-1)n - |W_j|+|V_i\cap W_j|$.

For $i\in [r]$, let $C(i):=\{j\in [t]: |V_i \cap W_j|\ne 0\}$ be the set of colors present in $V_i$.
For $j\in [t]$, let $\supp(j):=\{i\in [r]: |V_i \cap W_j|\ne 0\}$ be the set of parts that color $j$ is present.
Let
\begin{align*}
\Delta_j  := |W_j| - \min_{i\in \supp(j)} |V_i \cap W_j|.
%= \sum_{i=1}^r w_{i, j} - \min_{i\in \supp(j)} w_{i, j}.
%\Delta_j  := |W_j| - \min_{j\in \supp(j)\{ |V_i \cap W_j| : j \in [r] \text{ and } V_j \cap W_i \ne \emptyset \}.
\end{align*}
Thus $\delta(G)\leq \min_{j\in [t]} \{(r-1)n - \Delta_j \}$.
Hence we have for all $j \in [t]$,
\begin{align}
	\Delta_j < \Delta^* = \frac{(m-1)(r-1)n}{mt-2}.
	\label{eqn:Delta_i}
\end{align}

\begin{claim} \label{clm:Wj}
For all $j \in [t]$, $|W_j| < m \Delta^* /(m-1)$.
\end{claim}

\begin{proofclaim}
Note that 
\begin{align*}
\frac{m \Delta^*}{m-1} =  \frac{m (r-1)n}{mt-2} \geq \frac{m (m-1)t n }{mt-2}
 > (m-1)n. 
\end{align*}
If $|\supp(j)| \le m-1$, then $|W_j|\le (m-1)n < \frac{m}{m-1} \Delta^*$ as desired.
Hence we may assume that $|\supp(j)| \ge m$.
Using the definition of $\D_j$ and \eqref{eqn:Delta_i}, we obtain that
\begin{align*}
	\frac{m-1}{m} |W_j| \le \frac{|\supp(j)|-1}{|\supp(j)|} |W_j| \le \Delta_j < \Delta^*.
\end{align*}
Hence the claim follows. 
\end{proofclaim}

%For all $j \in [t]$, we know that $(r-1)n - \D^* < \delta(G) \le (r-1)n - \Delta_j$ and so
%\begin{align}
	%\Delta_j < \Delta^* = \frac{(m-1)(r-1)n}{r}< (m-1) n.
	%\label{eqn:Delta_i}
%\end{align}
%%%%%%%%%%%%%%%%%%%%%%%%%%%%%%%%%%%%%%%%%%%%%%%%%%%%%%%%%%%%%%%%%%%%%%%%%%%%%%%%%%%%%%%%%%%%%%%%%%%%%%%

Suppose that $|C(i)| =1$ for all $i \in [r]$. 
Every $V_i$ is a subset of some $W_j$ and consequently,  there exists $j \in [t]$ with $|W_j| \ge \lceil r/t \rceil n \ge m n$.
It follows that \[
 \Delta_j \ge |W_j|-n \ge (m-1)n \ge (m-1)\frac{ r-1}{mt-2}n = \Delta^*,
 \]
 a contradiction.

Without loss of generality, we assume that $|C(1)|=s\ge 2$.
For every $j \in C(1)$, we know that $|W_j| \le \Delta_j + |V_1 \cap W_j|$ from the definition of $\D_j$.
Hence, 
\begin{align*}
\sum_{j \in C(1)} |W_j| \le \sum_{j \in C(1)} ( \Delta_j + |V_1 \cap W_j|) \stackrel{(\ref{eqn:Delta_i})}{<} s \Delta^*  + n.
\end{align*}
Together with Claim~\ref{clm:Wj}, this gives 
\begin{align*}
	rn = \sum_{j \in [t]} |W_j| & < s \Delta^* + n + (t-s) \frac{m}{m-1} \Delta^* = \frac{mt-s}{m-1}\D^* + n \le  \frac{mt-2}{m-1}\D^* + n  = r n,
\end{align*}
a contradiction. 
\end{proof}

\section{Proof of the main results}\label{sec:main}

The proofs of Theorems~\ref{thm:-1} and~\ref{thm:-2} and Corollary~\ref{thm:t4} now follow easily from our auxiliary results.

\begin{proof}[Proof of Theorem~\ref{thm:-1}]
The $m=2$ case of the theorem is precisely Lemma~\ref{lma:f(n,2t-1,t+1)}.
For $m\geq 3$, we may apply
Corollary~\ref{cor:AES1} (with $a:=1$) to conclude that $f(n,r,t+1)=\delta(n,r,t)$. Then  Proposition~\ref{prop:d(n,r,t)3} implies $\delta(n,r,t)\leq (r-m)n=(r-\lceil r/t \rceil)n$; together with the lower bound in (\ref{eq:TuranBd}) this completes the proof.
\end{proof}

\begin{proof}[Proof of Theorem~\ref{thm:-2}]
Under Condition (i), we first apply Corollary~\ref{cor:AES1} to obtain that 
$f(n,r,t+1)=\delta(n,r,t)$. 
Then Propositions~\ref{prop:d(n,r,t)3lower} and~\ref{prop:d(n,r,t)3} give the desired lower and upper bounds, respectively.
Under Condition (ii), we apply Corollary~\ref{newcor} instead of Corollary~\ref{cor:AES1}.\COMMENT{\bf YZ expanded the proof. AL Rephrased the second sentence.}
%Corollary~\ref{newcor}, Proposition~\ref{prop:d(n,r,t)3lower} and Proposition~\ref{prop:d(n,r,t)3} together immediately imply Theorem~\ref{thm:-2}(ii).
%Theorem~\ref{thm:-2}(i) follows by replacing 
\end{proof}

\begin{proof}[Proof of Corollary~\ref{thm:t4}]
The case when  $r \equiv 0 \pmod{3}$ follows from \eqref{eq:BES}.
The case when $r \equiv 2 \pmod{3}$ follows immediately from Theorem~\ref{thm:-1}.
If $r \equiv 1 \pmod{3}$ then $r=3m-2$ for some $m \geq 4$. Thus Condition~(i) of Theorem~\ref{thm:-2} holds provided that $m \geq 6$. If $r=10$ and $n \geq 60$ or $r=13$ and $n \geq 22$, then it is easy to check that Condition~(ii) of  Theorem~\ref{thm:-2} holds. Thus Theorem~\ref{thm:-2} yields the corollary in this case.\COMMENT{YZ used Condition (i) or (ii) of Theorem~\ref{thm:-2} instead of Theorem~\ref{thm:-2}(i) or Theorem~\ref{thm:-2}(ii).}
\end{proof}

\section{Concluding remarks}\label{sec:conclude}
In this paper we have resolved Problem~\ref{pro:BES} for many choices of $r$ and $t$.
For the remaining open cases, it would be interesting to establish when (if at all) a lower bound construction from Section~\ref{sec:construct} is extremal.
One obvious case would be to determine $f(n,7,4)$, which is the only remaining case for $f(n,r,4)$. 

Our results show that $f(n,r,t+1)=\delta(n,r,t)$ when $r$ is large compared to $t$. It would be interesting to determine all values of $r$ and $t$ for which this equality holds. Proposition~\ref{prop:d(n,r,t)3} and \eqref{eq:rr} together show that $f(n, t+1, t+1)> \delta(n, t+1, t)$ when $t\ge 3$ is odd and $n$ is sufficiently large. Indeed, Proposition~\ref{prop:d(n,r,t)3} implies that $\delta(n, t+1, t)\le tn - \lceil \frac{tn}{2t-2}\rceil$ for every $t\ge 2$. If $t$ is odd, then by \eqref{eq:rr}, we have
\[
f(n, t+1, t+1)= tn - \left\lceil \frac{(t+1)n}{2t}\right\rceil > tn - \left\lceil \frac{tn}{2t-2}\right\rceil \ge \delta(n, t+1, t).
\]

As mentioned in the Introduction, Bollob\'as, Erd\H{o}s, and Straus~\cite{BEST} determined the largest $\delta(G)$ among all $K_3$-free (not necessarily balanced) $r$-partite graphs $G$ for all $r$. It is natural to extend the results in the present paper to \emph{unbalanced} multipartite graphs as well. 

It is also natural to ask for the largest $\delta(G)$ among $H$-free multipartite graphs $G$ for a fixed graph $H\ne K_t$. For example, Bollob\'as, Erd\H os, and Szemer\'edi~\cite{BES75} showed that if $G$ is a tripartite graph with $n$ vertices in each part and with $\delta(G)\ge n + \frac{1}{\sqrt{2}} n^{3/4}$, then $G$ contains a copy of $K_3(2)$; they asked if $\delta(G)\ge n + C n^{1/2}$ suffices. Furthermore, extending the aforementioned multipartite Tur\'an theorem of Bollob\'as, Erd\H{o}s, and Straus~\cite{BEST}, there has been recent work on determining the largest $e(G)$ among all multipartite graphs $G$ on $n$ vertices that contain no multiple (disjoint) copies of $K_t$, see, e.g., \cite{MR3952133, HZ19}.

\medskip
The following result might be useful for constructing extremal examples for the remaining open cases of Problem~\ref{pro:BES}. It shows that given an upper bound on $\Delta (n_0,r_0,t_0)/n_0$ one can obtain an upper bound on $\Delta (n,r,t)/n$ for  other triples $(n,r,t)$.
\begin{proposition}\label{generalconstruct}
Let $r_0, t_0 \in \mathbb N$ so that $2\leq t_0 \leq r_0 $. Let $\Delta_0\geq 0$ be such that $\Delta (n_0,r_0,t_0)/n_0 \leq \Delta _0$ for all $n_0 \in \mathbb N$.
Let $n,k \geq 2$ be integers. Set $r:=r_0k$ and $t:=k+t_0$.  
Then there exists an $r$-partite graph $G$ with parts of size $n$
%$V_1, \dots , V_{r}$ where $|V_1|= \cdots = |V_{r}|= n$ 
so that:
 \begin{itemize}
     \item $G$ contains no crossing independent set of size $t$; \COMMENT{YZ added an AND -- feel free to remove. AT: I think after a semi-colon you are not `suppose' to have an `and'}
     \item $\Delta (G) \leq (r_0-1)\left \lceil \frac{\Delta _0+(k-1)r_0}{\Delta_0+kr_0-1}  \cdot n \right \rceil  $.
 \end{itemize}
 That is,
 $$\Delta (n,r,t)\leq (r_0-1)\left \lceil \frac{\Delta _0+(k-1)r_0}{\Delta_0+kr_0-1} \cdot n \right \rceil ,$$
 or equivalently
 $$f(n,r,t)\geq (r-1)n- (r_0-1)\left \lceil \frac{\Delta _0+(k-1)r_0}{\Delta_0+kr_0-1} \cdot n \right \rceil.$$
\end{proposition}
%Note that, by (\ref{eq:Turan2}), in Proposition~\ref{generalconstruct}  one can set 
%$\Delta_0:= \lceil \frac{r_0}{t_0-1} \rceil -1$ if $t_0 \leq r_0$. 
\begin{proof}
\al{\bf defined $\ell$ to replace $\lfloor \alpha n \rfloor$ and other places in the proof}Let $\ell := \lfloor (r_0-1)n /(\Delta _0+kr_0-1) \rfloor$.
We construct an $r$-partite graph $G$ with parts $V_1, \dots , V_{r}$ of size $n$ as follows.\COMMENT{YZ shortened the sentence.}
%where $|V_1|= \cdots = |V_{r}|= n$ and that contains no crossing independent set of size $t$.
Partition each $V_i$ into $L_i\cup S_i$ such that $|S_i| = \ell$ and $|L_i| = n-\ell$. We call the vertices in each $L_i$ \emph{large} and those in each $S_i$ \emph{small}.
We partition $[r]$ into $k$ blocks of size $r_0$ by assigning $i$ and $j$ to the same block if $\lceil i/r_0 \rceil = \lceil j/r_0 \rceil$. 
We form a complete bipartite graph between $L_i$ and $L_j$ (\emph{joining $L_i$ and $L_j$} for short) if and only if $i$ and $j$ are in the same block. We join $S_i$ and $S_j$  if $i$ and $j$ are not in the same block. 
\al{\bf rephrased the following sentence.}
By the definition of $\Delta_0$, there exists an $r_0$-partite graph $G_S$ with $\ell$ vertices in each part and $\Delta(G_S)  \leq  \Delta _0 \ell$ containing no crossing independent set of size~$t_0$.
%Additionally, let $G_S$ be an $r_0$-partite graph with
%$\ell$ vertices in each part 
%with $\Delta (G_S) \leq  \Delta _0 \ell$ with no crossing independent set of size~$t_0$;\footnote{YZ: this sentence does not read well AL shorten this sentence a bit.} such a graph $G_S$ exists by definition of $\Delta _0$. 
We place a copy of $G_S$ in each block so that the sets $S_i$ in the block each form one of the vertex classes of~$G_S$.

We claim that $G$ contains no crossing independent set of size $t$. Indeed,  a crossing independent set contains at most one large vertex from each block, and at most $t_0-1$  small vertices. Thus the largest crossing independent set has size $k+t_0-1 = t-1$.\al{\bf was $<t$.}

Note that 
$$\Delta (G)= \max \{ (\Delta _0+(k-1)r_0) \ell, (r_0-1) (n-\ell).$$
The choice of $\ell$ ensures
$$\Delta (G)\leq (r_0-1)\left \lceil \frac{\Delta _0+(k-1)r_0}{\Delta_0+kr_0-1} \cdot n \right \rceil ,$$
as desired.
\end{proof}


\begin{thebibliography}{10}

\bibitem{MR955135}
N.~Alon.
\newblock The linear arboricity of graphs.
\newblock {\em Israel J. Math.}, 62(3):311--325, 1988.

\bibitem{MR0340075}
B.~Andr\'{a}sfai, P.~Erd\H{o}s, and V.~T. S\'{o}s.
\newblock On the connection between chromatic number, maximal clique and
  minimal degree of a graph.
\newblock {\em Discrete Math.}, 8:205--218, 1974.

\bibitem{MR3952133}
P.~Bennett, S.~English, and M.~Talanda-Fisher.
\newblock Weighted {T}ur\'{a}n problems with applications.
\newblock {\em Discrete Math.}, 342(8):2165--2172, 2019.

\bibitem{BEST}
B.~Bollob\'as, P.~Erd\H{o}s, and E.~Straus.
\newblock Complete subgraphs of chromatic graphs and hypergraphs.
\newblock {\em Utilitas Math.}, 6:343--347, 1974.

\bibitem{BES75}
B.~Bollob\'{a}s, P.~Erd\H{o}s, and E.~Szemer\'{e}di.
\newblock On complete subgraphs of {$r$}-chromatic graphs.
\newblock {\em Discrete Math.}, 13(2):97--107, 1975.

\bibitem{bondy}
A.~Bondy, J.~Shen, S.~Thomass\'e, and C.~Thomassen.
\newblock Density conditions for triangles in multipartite graphs.
\newblock {\em Combinatorica}, 26:121--131, 2006.

\bibitem{MR0329905}
P.~Erd\H{o}s.
\newblock Unsolved problems.
\newblock In D.~J.~A. Welsh and D.~R. Woodall, editors, {\em Proceedings of the
  {C}onference on {C}ombinatorial {M}athematics held at the {M}athematical
  {I}nstitute, {O}xford, 3--7 {J}uly 1972}, pages 351--363. The Institute of
  Mathematics and its Applications, Southend-on-Sea, 1972.

\bibitem{HZ19}
J.~{Han} and Y.~{Zhao}.
\newblock {Tur\'an number of disjoint triangles in 4-partite graphs}.
\newblock {\em arXiv e-prints}, page arXiv:1906.01812, Jun 2019.

\bibitem{MR2195582}
P.~Haxell and T.~Szab\'{o}.
\newblock Odd independent transversals are odd.
\newblock {\em Combin. Probab. Comput.}, 15(1-2):193--211, 2006.

\bibitem{MR1860440}
P.~E. Haxell.
\newblock A note on vertex list colouring.
\newblock {\em Combin. Probab. Comput.}, 10(4):345--347, 2001.

\bibitem{MR1208805}
G.~P. Jin.
\newblock Complete subgraphs of {$r$}-partite graphs.
\newblock {\em Combin. Probab. Comput.}, 1(3):241--250, 1992.

\bibitem{pfender}
F.~Pfender.
\newblock Complete subgraphs in multipartite graphs.
\newblock {\em Combinatorica}, 32:483--495, 2012.

\bibitem{MR2246152}
T.~Szab\'{o} and G.~Tardos.
\newblock Extremal problems for transversals in graphs with bounded degree.
\newblock {\em Combinatorica}, 26(3):333--351, 2006.

\bibitem{Tu41}
P.~Tur\'an.
\newblock On an extremal problem in graph theory.
\newblock {\em Mat. Fiz. Lapok}, 48:436--452, 1941.

\end{thebibliography}
\end{document}